\documentclass[fleqn]{mat01}
\usepackage{times,mathtimy,amssymb,latexsym,boxit}
\begin{document}

\setcounter{page}{431}
\firstpage{431}

\font\xxx=msam10
\def\ab{\mbox{\xxx{\char'245}}}

\newcommand{\R}{\mathbb{R}}
\newcommand{\CA}{{\mathcal{A}}}

\newtheorem{theo}{Theorem}
\renewcommand\thetheo{\arabic{section}.\arabic{theo}}
\newtheorem{theor}[theo]{\bf Theorem}
\newtheorem{lem}[theo]{Lemma}
\newtheorem{propo}{\rm PROPOSITION}
\newtheorem{rmk}[theo]{Remark}
\newtheorem{defn}[theo]{\rm DEFINITION}
\newtheorem{exam}{Example}
\newtheorem{coro}[theo]{\rm COROLLARY}
\def\conjecture{\trivlist\item[\hskip\labelsep{\it Conjecture.}]}

\renewcommand{\theequation}{\thesection\arabic{equation}}

\title{A probabilistic approach to second order variational inequalities
with bilateral constraints}

\markboth{Mrinal K Ghosh and K~S Mallikarjuna Rao}{A probabilistic approach to second order variational
inequalities}

\author{MRINAL K GHOSH$^{*,\dagger}$ and K~S MALLIKARJUNA RAO$^{\ddagger}$}

\address{$^{*}$Department of Mathematics, Indian Institute of Science,
Bangalore~560~012, India\\
\noindent $^{\dagger}$Department of Electrical and Computer Engineering,
University of Texas, Austin, TX~78712, USA\\
\noindent $^{\ddagger}$CMI, Universit\'{e} de Provence, 39, Rue F. J.
Curie, 13 453 Marseille, France\\
\noindent Email: mkg@math.iisc.ernet.in; mrinal@ece.utexas.edu}

\volume{113}

\mon{November}

\parts{4}

\Date{MS received 5 April 2002; revised 8 May 2003}

\begin{abstract}
We study a class of second order variational inequalities with bilateral
constraints. Under certain conditions we show the existence of a unique
viscosity solution of these variational inequalities and give a
stochastic representation to this solution. As an application, we study
a stochastic game with stopping times and show the existence of a saddle
point equilibrium.
\end{abstract}


\maketitle

\section{Introduction and preliminaries}

We study a class of second order nonlinear variational
inequalities with bilateral constraints. This type of inequalities
arises in zero sum stochastic differential games of mixed type
where each player uses both continuous control and stopping times.
Under a non-degeneracy assumption Bensoussan and Friedman
\cite{BL,Fr}  have studied this type of  problems. They proved the existence
of a unique solution of these variational inequalities in certain
weighted Sobolev spaces. This result together with certain
techniques from stochastic calculus is then applied to show that
the unique solution of these inequalities is the value function of
certain stochastic differential games of mixed type. In this
paper we study the same class of variational inequalities without
the non-degeneracy assumption. The non-degeneracy assumption is
crucially used in the analysis of the problem in \cite{BL,Fr}.
Thus the method used in \cite{BL,Fr} does not apply to the
degenerate case. We study the problem via the theory of viscosity
solutions. We transform the variational inequalities with
bilateral constraints to Hamilton--Jacobi--Isaacs (HJI for short)
equations associated with a stochastic differential game problem with
continuous control only. Then using standard results from the
theory of viscosity solutions, we show that the value function of
this stochastic differential game with continuous control is the
unique viscosity solution of the corresponding variational
inequalities. Then for a special case  we identify this unique
viscosity solution as the value function of the stochastic game
with stopping times. We now describe our problem.

Let $U_i, i=1,2$, be the compact metric spaces.  Let
\setcounter{equation}{0}
\begin{equation*}
b: \R^d \times U_1 \times U_2 \to \R^d
\end{equation*}
and
\begin{equation*}
a: \R^d \times U_1 \times U_2 \to \R^{d \times d}.
\end{equation*}
We assume that:\vspace{.6pc}

\noindent (A1)\vspace{.6pc}

\noindent The functions $b$ and $a$ are
bounded and continuous, $a(x, u_1, u_2)$ is $C^2$ in $x$ uniformly
with respect to $u_1, u_2$.
The matrix
$a$ is symmetric and non-negative definite. Further there exists
constant
$C_1 >0$ such that for all $u_i \in U_i, i = 1, 2$,
\begin{equation*}
|b(x,u_1,u_2) - b(y, u_1, u_2) |  \leq C_1 |x-y|.
\end{equation*}
Let
\begin{equation*}
r:\R^d \times U_1 \times U_2 \to \R
\end{equation*}
and
\begin{equation*}
\psi_i : \R^d \to \R, \,\, i=1,2.
\end{equation*}
We assume that\vspace{.6pc}

\noindent (A2)
\begin{enumerate}
\renewcommand{\labelenumi}{(\roman{enumi})}
\leftskip .5pc
\item $r, \psi_1, \psi_2$ are bounded and continuous.

\item There exists a  constant $C_2 > 0$ such
that for all $x,y \in \R^d, (u_1, u_2) \in U_1 \times U_2$,
\begin{align*}
&|r(x,u_1,u_2) - r(y,u_1,u_2)| + |\psi_1(x) - \psi_1(y)|
+|\psi_2(x) - \psi_2(y)|\\
&\quad\ \leq C_2 |x-y|.
\end{align*}
\item $\psi_2 \leq \psi_1$.
\end{enumerate}

\noindent Let $H^{+}, H^{-} : \R^d \times \R^d \times \R^{d\times d} \to \R$
be defined by
\begin{align*}
H^+(x, p, X) &= \inf_{u_1 \in U_1} \sup_{u_2 \in U_2}
\left[ \frac{1}{2}{\rm tr}(a(x,u_1,u_2)X) \right.\\
&\quad\  +  b(x,u_1,u_2) \cdot p + r(x,u_1,u_2)\bigg\rbrack, \\
H^-(x, p, X) &= \sup_{u_2 \in U_2} \inf_{u_1 \in U_1}
\left[ \frac{1}{2}{\rm tr}(a(x,u_1,u_2)X) \right.\\
&\quad\  + b(x,u_1,u_2) \cdot p + r(x,u_1,u_2)\bigg\rbrack.
\end{align*}
Consider the following Hamilton--Jacobi--Isaacs variational
inequalities with bilateral constraints
\begin{equation}\label{St2s1}
\left. \begin{array}{ll}
\psi_2(x) \leq v(x) \leq \psi_1(x), & \forall ~  x\\[.2pc]
\lambda v(x) - H^+(x,Dv(x),D^2 v(x)) = 0, & \mbox{if } \psi_2(x) < v(x) < \psi_1(x) \\[.2pc]
\lambda v(x) - H^+(x,Dv(x), D^2 v(x)) \geq 0, & \mbox{if } v(x) = \psi_2(x) \\[.2pc]
\lambda v(x) - H^+(x,Dv(x), D^2 v(x)) \leq 0, & \mbox{if } v(x) = \psi_1(x) \end{array} \right\}
\end{equation}
and
\begin{equation} \label{St2s2}
\left. \begin{array}{ll}
\psi_2(x) \leq v(x) \leq \psi_1(x), & \forall  ~x\\[.2pc]
\lambda v(x) - H^-(x,Dv(x), D^2 v(x)) = 0, & \mbox{if }
\psi_2(x) < v(x) < \psi_1(x) \\[.2pc]
\lambda v(x) - H^-(x,Dv(x), D^2 v(x)) \geq 0, & \mbox{if }
v(x) = \psi_2(x) \\[.2pc]
\lambda v(x) - H^-(x,Dv(x), D^2 v(x)) \leq 0, & \mbox{if } v(x) =
\psi_1(x)
\end{array} \right\}.
\end{equation}

By a classical solution of (\ref{St2s1}), we mean a $C^2$-function
$v$ satisfying (\ref{St2s1}). Similarly a classical solution of
(\ref{St2s2}) is defined.

The rest of our paper is structured as follows. In \S2, we
introduce the notion of viscosity solution and  establish the
existence of  unique viscosity solutions of these variational
inequalities by a probabilistic method. In \S3, we apply
these variational inequalities  to treat a stochastic game with
stopping times. We establish the existence of a saddle point
equilibrium for this problem. Section 5 contains some concluding
remarks.

\section{Viscosity solutions}

To motivate the definition of viscosity solutions of the
variational inequalities we first prove the following result.

\setcounter{equation}{0}

\begin{theor}[\!]\label{St2t1}
Assume {\rm (A2)(iii)}. A function $v \in C^2(\R^d)$ is a classical solution
of {\rm (\ref{St2s1})} if and only if it is a classical
solution of the equation
\begin{align}\label{St2s3}
&\max \{ \min \{ \lambda v(x) - H^+(x, Dv(x),D^2 v(x)); \lambda(v(x) -
\psi_2(x)) \};\nonumber\\
&\hskip 1cm \lambda(v(x) - \psi_1(x)) \} = 0.
\end{align}
Similarly a function $v \in C^2(\R^d)$ is a classical
solution of {\rm (\ref{St2s2})} if\, and\, only\, if\, it\, is a classical
solution of the equation
\begin{align}\label{St2s4}
&\min \{ \max \{ \lambda v(x) - H^-(x, Dv(x),D^2 v(x)); \lambda(v(x)
- \psi_1(x)) \};\nonumber\\
&\hskip 1cm \lambda (v(x) - \psi_2(x)) \} = 0.
\end{align}
\end{theor}

\begin{proof}
Let $v$ be a classical solution of (\ref{St2s1}).  Suppose
$x$ is such that $\psi_2(x) < v(x) < \psi_1(x)$.  Then
\begin{align*}
&\lambda v(x) - H^+(x, Dv(x), D^2 v(x)) = 0,~~ v(x) - \psi_2(x) > 0,\\
&\hskip 1cm v(x) - \psi_1(x) < 0.
\end{align*}
Thus  (\ref{St2s3}) clearly holds in this case.  Now
if $v(x) = \psi_2(x)$, then
\begin{equation*}
\min \{ \lambda v(x) - H^+(x, Dv(x), D^2 v(x)) ; \lambda (v(x) - \psi_2(x))
\} = 0,
\end{equation*}
and hence (\ref{St2s3}) is satisfied.  Finally assume
$v(x) = \psi_1(x)$, then
\begin{equation*}
\min \{ \lambda v(x) - H^+(x, Dv(x), D^2 v(x)) ; \lambda (v(x) - \psi_2(x))
\} \leq 0,
\end{equation*}
and hence
\begin{align*}
&\max\{ \min \{ \lambda v(x) - H^+(x, Dv(x), D^2 v(x)) ; \lambda (v(x)
- \psi_2(x)) \};\\
&\hskip 1cm \lambda (v(x) - \psi_1(x)) \} = 0.
\end{align*}
Thus $v$ satisfies (\ref{St2s3}).
We now show the converse.  It is clear from (\ref{St2s3})
that $v(x) \leq \psi_1(x)$.  If $v(x) = \psi_1(x)$ for
some $x$, then it clearly satisfies $v(x) \geq \psi_2(x)$
by (A2)(iii).
Now let $v(x) < \psi_1(x)$.  Then from (\ref{St2s3}), we
have
\begin{equation*}
\min \{ \lambda v(x) - H^+(x, Dv(x), D^2 v(x)) ; \lambda (v(x) - \psi_2(x))
\} = 0
\end{equation*}
and hence $v(x) - \psi_2(x) \geq 0$.  Thus for all $x$,
we have $\psi_2(x) \leq v(x) \leq \psi_1(x)$.  Now let
$v(x) <\psi_1(x)$.  Then from the above equation, we have
\begin{equation*}
\lambda v(x) - H^+(x, Dv(x), D^2 v(x)) \geq 0.
\end{equation*}
Similarly if $v(x) > \psi_2(x)$, we can show that
\begin{equation*}
\lambda v(x) - H^+(x, Dv(x), D^2 v(x)) \leq 0.
\end{equation*}
Thus $v$ is a classical solution of (\ref{St2s1}). This concludes the
proof of the first part. The second part of the theorem can be proved in
a similar way.\hfill \ab
\end{proof}

\begin{rmk}
{\rm Under (A2)(iii), we can also show that a function $v \in
C^2(\R^d)$ is a classical solution of (\ref{St2s1}) if and only if
it is a classical solution of the equation
\begin{align*}
&\min \{ \max \{ \lambda v(x) - H^+(x, Dv(x), D^2 v(x)); \lambda(v(x) -
\psi_1(x)) \};\\
&\hskip 1cm \lambda (v(x) - \psi_2(x)) \} = 0.
\end{align*}
Similarly a function $v \in C^2(\R^d)$ is a classical solution of
(\ref{St2s2}) if and only if it is a classical solution of the
equation
\begin{align*}
&\max \{ \min \{ \lambda v(x) - H^-(x, Dv(x), D^2 v(x)); \lambda(v(x) -
\psi_2(x)) \};\\
&\hskip 1cm \lambda (v(x) - \psi_1(x)) \} = 0.
\end{align*}}
\end{rmk}

Theorem \ref{St2t1} motivates us to define viscosity solutions for
(\ref{St2s1}) and (\ref{St2s2}) using (\ref{St2s3}) and (\ref{St2s4})
repsectively.

\begin{defn}$\left.\right.$\vspace{.5pc}

\noindent {\rm An upper semicontinuous function $v: \R^d \to \R$
is said to be a viscosity subsolution of (\ref{St2s1})
if it is a viscosity subsolution of (\ref{St2s3}).
Similarly a lower semicontinuous function $v: \R^d \to \R$
is said to be a viscosity supersolution of (\ref{St2s1})
if it is a viscosity supersolution of (\ref{St2s3}).
A function which is both sub- and super-solution of
(\ref{St2s1}) is called a viscosity solution of (\ref{St2s1}).
Similarly, viscosity sub- and super-solutions
of (\ref{St2s2}) are\break defined.}
\end{defn}

We now address the  question of showing the existence of unique
viscosity solutions of (\ref{St2s1}) and (\ref{St2s2}). This is
done by showing that (\ref{St2s1}) and (\ref{St2s2}) are
equivalent to Hamilton--Jacobi--Isaacs equations corresponding to a
stochastic differential game.

Let $\omega_1, \omega_2$ be two symbols. We formulate a zero sum
stochastic differential. In this game, $\bar{U}_i$ is the set of
controls for player $i$, where $\bar{U}_i = U_i \cup \{\omega_i\},
i = 1,2$. Let $\sigma(\cdot, \cdot, \cdot)$ be the non-negative
square root of $a(\cdot, \cdot, \cdot).$
Extend $b,\sigma, r$ to
\begin{align*}
&\bar{b} : \R^d \times \bar{U}_1 \times \bar{U}_2 \to \R^d,~~
\bar{\sigma} : \R^d \times \bar{U}_1  \times
\bar{U}_2 \to \R^{d \times d},\\
&\bar{r} :\R^d \times  \bar{U}_1 \times \bar{U}_2 \to \R,
\end{align*}
respectively such that
\begin{align*}
&\bar{b}(x,\omega_1, \cdot) \equiv 0,~~ \bar{b}(x, \cdot, \omega_2)
\equiv 0,~~ \bar{\sigma}(x,\omega_1, \cdot) \equiv 0, ~~ \bar{\sigma}
(x, \cdot, \omega_2)\equiv 0,\\
&\bar{r}(x, \omega_1, u_2) = \lambda \psi_1(x) \mbox{
for all } u_2 \in U_2
\mbox{ and } \bar{r}(x, \cdot, \omega_2)
\equiv \lambda \psi_2(x).
\end{align*}
Let $({\Omega}, {{\mathcal F}}, {P})$ be a complete probability
space and ${W}(t)$ a standard $d$-dimensional Brownian motion
on it. Let $\bar{\CA}_i$ denote the set of all $\bar{U}_i$-valued
functions progressively measurable with respect to the
process ${W}(t)$.  For $(\bar{u}_1(\cdot),\bar{u}_2(\cdot))
\in \bar{\CA}_1 \times \bar{\CA}_2$, consider the controlled
stochastic differential equation
\begin{equation}\label{St2s5}
\left. \begin{array}{rcl}
{\rm d}\bar{X}(t) & = &
\bar{b}(\bar{X}(t), \bar{u}_1(t),\bar{u}_2(t)) \,\, {\rm d}t
+ \bar{\sigma} (\bar{X}(t), u_1(t), u_2(t)) \,\, {\rm d}{W}(t) \\[.2pc]
\bar{X}(0) & = &x
\end{array}
\right\rbrace.
\end{equation}
Let the payoff function be defined by
\begin{equation*}
\bar{R}(x, \bar{u}_1(\cdot), \bar{u}_2(\cdot)) = E \left[ \int_0^\infty
{\rm e}^{- \lambda t} \bar{r} (\bar{X}(t), \bar{u}_1(t), \bar{u}_2(t))
\,\, {\rm d}t \right].
\end{equation*}
A strategy for the player $1$ is a `non-antiticipating' map
$\alpha : \bar{\CA}_2 \to \bar{\CA}_1$, i.e.,  for any $u_2,
\bar{u}_2 \in \CA_2$ such that $u_2(s) = \bar{u}_2(s)$ for all $0
\leq s \leq t$ then we have $\alpha[u_2](s) =
\alpha[\bar{u}_2](s), 0 \leq s \leq t$. Let $\bar{\Gamma}$ denote
the set of all non-anticipating strategies for player $1$. Similarly
strategies for player $2$ are defined. Let the set of all
non-anticipating strategies for player $2$ be denoted by
$\bar{\Delta}$. Then the upper and lower value functions are
defined by
\begin{align*}
\bar{V}^+(x) &= \sup_{\beta \in \bar{\Delta} }
\inf_{\bar{u}_1\in \bar{\CA}_1} \bar{R}(x, \bar{u}_1(\cdot),
\beta[\bar{u}_1](\cdot)),\\
\bar{V}^-(x) &= \inf_{\alpha \in \bar{\Gamma} }
\sup_{\bar{u}_2\in \bar{\CA}_2} \bar{R}(x, \alpha[\bar{u}_2](\cdot),
\bar{u}_2(\cdot)).
\end{align*}
Then we can closely follow the arguments in \cite{FS} to
show that under (A1) and (A2),
$\bar{V}^+$ and $\bar{V}^-$, respectively, are  unique viscosity
solutions of
\begin{equation}\label{St2s6}
\lambda v(x) - \bar{H}^+(x, Dv(x), D^2 v(x)) = 0
\end{equation}
and
\begin{equation}\label{St2s7}
\lambda v(x) - \bar{H}^-(x, Dv(x), D^2 v(x)) = 0
\end{equation}
in the class of bounded continuous functions, where $\bar{H}^+,
\bar{H}^- : \R^d \times \R^d \times \R^{d \times d} \to \R$
are defined as follows:
\begin{align*}
\bar{H}^+(x, p, X) &= \inf_{\bar{u}_1 \in \bar{U}_1} \sup_{\bar{u}_2 \in
\bar{U}_2} \left[ \frac{1}{2}{\rm tr}(\hat{a}(x, u_1, u_2)X) \right.\\
&\quad\ + \bar{b}(x, \bar{u}_1, \bar{u}_2) \cdot p + r(x,
\bar{u}_1,\bar{u}_2) \bigg\rbrack, \\ \bar{H}^-(x, p, X) &=
\sup_{\bar{u}_2 \in \bar{U}_2} \sup_{\bar{u}_1 \in \bar{U}_1}
\left[\frac{1}{2} {\rm tr}(\hat{a}(x, u_1, u_2)X) \right.\\
&\quad\ + \bar{b}(x, \bar{u}_1, \bar{u}_2) \cdot p + r(x,
\bar{u}_1,\bar{u}_2) \bigg\rbrack,
\end{align*}
where $\hat{a} = \bar{\sigma} \bar{\sigma}^*$.
We now establish the equivalence of (\ref{St2s1}) and
(\ref{St2s2}) with (\ref{St2s6}) and (\ref{St2s7})
respectively.

\begin{theor}[\!]\label{St2t5}
Assume {\rm (A2)(iii)}. A continuous function $v:\R^d \to \R$ is a viscosity
solution of {\rm (\ref{St2s6})} if and only if it is a
viscosity solution of {\rm (\ref{St2s1})}. Similarly
a continuous function $v:\R^d \to \R$ is a viscosity
solution of {\rm (\ref{St2s7})} if and only if it is a
viscosity solution of {\rm (\ref{St2s2})}.
\end{theor}

\begin{proof}
The proof of this theorem is a simple consequence
of the observation that
\begin{equation}\label{St2s8}
\bar{H}^+(x,p, X) = (H^+(x,p, X) \vee \lambda \psi_2(x) )
\wedge \lambda \psi_1(x)
\end{equation}
and
\begin{equation}\label{St2s9}
\bar{H}^-(x,p, X) = (H^-(x,p, X) \wedge \lambda \psi_1(x) )
\vee \lambda \psi_2(x).
\end{equation}
\hfill \ab
\end{proof}

As a consequence of this  result we have the following
existence and uniqueness result for the solutions of
(\ref{St2s1}) and (\ref{St2s2}).

\begin{coro}\label{St2t6}$\left.\right.$\vspace{.5pc}

\noindent Assume {\rm (A1)} and {\rm (A2)}. Then $\bar{V}^+$ and $\bar{V}^-$ are unique
vicosity solutions of {\rm (\ref{St2s1})} and {\rm (\ref{St2s2})} respectively
in the class of bounded continuous functions.
\end{coro}

\begin{proof}
Since (\ref{St2s6}) has a unique viscosity solution in the
class of bounded continuous functions given by $V^+$, we get by
Theorem \ref{St2t5}, that $\bar{V}^+$ is the unique viscosity
solution of (\ref{St2s1}) in the class of bounded continuous
functions. Similarly $V^-$ is the unique viscosity solution of
(\ref{St2s2}) in the class of bounded continuous functions.\hfill \ab
\end{proof}

\begin{rmk}
{\rm In the classical case, it is quite clear from the proof of Theorem
\ref{St2t1} that under (A2)(iii), $\psi_2 \leq v \leq \psi_1$ if
$v$ is a classical solution of (\ref{St2s3}).  In fact this
remains true even in the case of viscosity solutions.
Indeed assume (A2)(iii). Then any  viscosity solution of
(\ref{St2s3}) satisfies
\begin{equation*}
\psi_2(x) \leq v(x) \leq \psi_1(x), \mbox{ for
all } x \in \R^d.
\end{equation*}
Similar result holds for  viscosity solutions of (\ref{St2s4}).
Observe that $\psi_2$ is a viscosity sub-solution of (\ref{St2s3})
and $\psi_1$ is a viscosity super-solution of (\ref{St2s3}). Thus the
desired result follows from a general comparison principle on
viscosity solutions \cite{CIL}.}
\end{rmk}

\section{Application to stochastic games}

\setcounter{theo}{0}
In this section we consider a stochastic  game with stopping
times.  We show the existence of a value and a saddle point
equilibrium for this problem.  We now describe the stochastic game
with stopping times.

Let $b : \R^d \to \R^d$ and $\sigma : \R^d \to \R^{d \times d}$.
Assume  that $b, \sigma$ are bounded and Lipschitz continuous.
Consider the stochastic differential equation
\setcounter{equation}{0}
\begin{equation} \label{St3s1}
\left. \begin{array}{rcl}
{\rm d}X(t) & = & b(X(t)) \,\, {\rm d}t + \sigma(X(t)) \,\, {\rm d}W(t),
\,\,t > 0\\[.2pc]
X(0) & =  & x
\end{array} \right\}.
\end{equation}
Here $W$ is a $d$-dimensional Brownian motion on an underlying
complete probability space $(\Omega, {\mathcal F}, P)$. Let $r,
\psi_1, \psi_2: \R^d \to \R$ be  bounded and  Lipschitz continuous
functions and $\psi_1 \geq \psi_2 $. Let $\lambda > 0$.  Define
\begin{align}\label{St3s2}
R(x, \theta, \tau) &= E \left[ \int_0^{\theta \wedge \tau} {\rm e}^{-\lambda s}
r(X(s))\,\, {\rm d}s + {\rm e}^{-\lambda (\theta \wedge \tau)}
[\psi_1(X(\theta)) \chi_{\theta < \tau} \right.\nonumber\\
&\quad\ + \psi_2(X(\tau)) \chi_{\tau \leq \theta}] \bigg\rbrack,
\end{align}
where $\theta, \tau$ are the stopping times with respect to the
$\sigma$-field generated by $W(t)$.

Let $\tilde{W}$ denote another $d$-dimensional Brownian motion
independent of $W$, constructed on an augmented probability space
which is also denoted by $(\Omega, {\mathcal F}, P)$ by an abuse of
notation.  Let $\sigma^\gamma : \R^d \to \R^{d \times 2d}$ be
defined by $ \sigma^{\gamma} = [ \sigma ~~ \gamma I_d]$, where
$I_d$ is the $d \times d$ identity matrix.  Now consider the
following stochastic differential equation
\begin{equation} \label{St3s3}
\left. \begin{array}{rcl}
{\rm d}X^{\gamma}(t) &  =  & b(X^{\gamma}(t)) \,\, {\rm d}t +
\sigma^\gamma(X^{\gamma}(t)) \,\,
{\rm d}\overline{W}(t), \,\,t > 0\\[.2pc]
X^{\gamma}(0) & =  & x
\end{array}
\right\rbrace,
\end{equation}
where $\overline{W} = [W, \tilde{W}]^*$. Player 1 tries to minimize
$R(x, \cdot, \cdot)$, as in (\ref{St3s2}), over stopping times $\theta$
(with respect to the $\sigma$-field generated by $\overline{W}(t)$),
whereas player 2 tries to maximize the same over stopping times $\tau$
(with respect to the $\sigma$-field generated by $\overline{W}(t)$).
Note that these stopping times need not be finite a.s.. In other words,
each player has the option of not stopping the game at any time. We
refer to \cite{St} for a similar treatment to stochastic games with
stopping times. We now define the lower and upper value functions. Let
\begin{align*}
V^-(x) = \sup_{\tau \geq 0} \inf_{\theta \geq 0} R(x, \theta, \tau),\\
V^+(x) = \inf_{\theta \geq 0} \sup_{\tau \geq 0} R(x, \theta, \tau).
\end{align*}
The stochastic game with stopping times is said to have a value if
$V^+ \equiv V^-$.

We now establish the existence of a value for this problem. Let $H
: \R^d \times \R^d \times \R^{d \times d} \to \R$ be defined by
\begin{equation*}
H(x,p, X) =  \frac{1}{2}\ {\rm tr}\ (a(x)X) + b(x) \cdot p + r(x),
\end{equation*}
where $a = \sigma \sigma^*$.

\begin{theor}[\!]
The stochastic game with stopping times has a value and a saddle
point equilibrium. The value of this game is the unique viscosity
solution in the class of bounded and continuous functions of the
variational inequalities with bilateral constraints given by
\begin{equation} \label{St3s4}
\left.
\begin{array}{ll}
\psi_2(x) \leq w(x) \leq \psi_1(x), & \forall x \\[.2pc]
\lambda w(x) - H(x, Dw(x), D^2 w(x)) = 0, &
\mbox{if } \psi_2(x) < w(x) < \psi_1(x) \\[.2pc]
\lambda w(x) - H(x, Dw(x), D^2 w(x)) \geq  0, &
\mbox{if }  w(x) = \psi_2(x)  \\[.2pc]
\lambda w(x) - H(x, Dw(x), D^2 w(x)) \leq  0, & \mbox{if }  w(x) =
\psi_1(x)
\end{array} \right\rbrace.
\end{equation}
\end{theor}

\begin{proof}
Using the results of \S2, it follows that there is a
unique viscosity solution $w$ of (\ref{St3s4}) in the class of
bounded continuous functions. We now identify this solution as the
value function of the stochastic game with stopping times.

We first assume that $w$ is $C^2$. Let $w(x) > \psi_2(x)$. Let
$\theta$ be any stopping times.  Define the stopping time
$\hat{\tau}$ by
\begin{equation*}
\hat{\tau} = \inf \{ t\geq 0 :\psi_2(X(t)) = w(X(t)) \},
\end{equation*}
where $X(t)$ is a solution of (\ref{St3s1}) with the initial
condition $X(0) = x$. Since $w$ is a smooth viscosity solution of
(\ref{St3s4}), by Ito's formula, we have  for any $T > 0$ and
stopping time $\theta$,
\begin{equation*} w(x) \leq E \left\lbrace \int_0^{T
\wedge \theta \wedge \hat{\tau}} {\rm e}^{-\lambda s} r(X(s)) \,\, {\rm d}s +
{\rm e}^{-\lambda (T \wedge \theta \wedge \hat{\tau})} w(X(T \wedge \theta
\wedge \hat{\tau}) \right\rbrace.
\end{equation*}
Letting $T \to \infty$ in the above equation, we obtain
\begin{equation*}
w(x) \leq E \left\lbrace \int_0^{ \theta \wedge \hat{\tau}} {\rm e}^{-\lambda
s} r(X(s)) \,\, {\rm d}s + {\rm e}^{-\lambda (\theta \wedge \hat{\tau})}
w(X(\theta \wedge \hat{\tau})) \right\rbrace.
\end{equation*}
Now using the first inequality in (\ref{St3s4}) and the definition of
$\hat{\tau}$ in the above, it follows that
\begin{equation*}
w(x) \: \leq \: R(x, \theta, \hat{\tau}).
\end{equation*}
Since $\theta$ is arbitrary, we get
\begin{equation*}
w(x) \leq V^-(x).
\end{equation*}

Next let $w$ is $C^2$ and $w(x) < \psi_1(x)$.  Define the
stopping time $\hat{\theta}$ by
\begin{equation*}
\hat{\theta} = \inf \{ t\geq 0 :\psi_1(X(t)) = w(X(t)) \},
\end{equation*}
where $X(t)$ is a solution of (\ref{St3s1}) with the initial
condition $X(0) = x$. Now using the foregoing arguments we can show
that
\begin{equation*}
w(x) \geq V^+(x).
\end{equation*}

Thus
\begin{equation*}
w \equiv V^+ \equiv V^-.
\end{equation*}

We now prove this result for the general case. Let $w_\epsilon$ be
the sup-convolution of $w$, i.e.,
\begin{equation*}
w_\epsilon (x) = \sup_{\xi \in \R^d} \left\lbrace w(\xi) - \frac{|\xi -
x|^2}{2 \epsilon} \right\rbrace.
\end{equation*}
Then $w_\epsilon \to w$ uniformly in $\R^d$ as
$\epsilon \to 0$, $w_\epsilon$ are bounded, Lipschitz
continuous, semiconvex and satisfy a.e. on $\R^d$,
\begin{equation*}
\lambda w_\epsilon(x) - H(x, D w_\epsilon(x),
D^2 w_\epsilon(x)) \leq \rho_0 (x),
\end{equation*}
for a modulus $\rho_0$ (see \cite{S2}). Now let
$w_\epsilon^\delta$ be the standard mollification of $w_\epsilon$.
Then $w_\epsilon^\delta$ are $C^{2,1}$, $w_\epsilon^\delta \to
w_\epsilon$ uniformly in $\R^d$ and $Dw_\epsilon^\delta(x) \to
Dw_\epsilon(x), D^2w_\epsilon^\delta (x) \to D^2 w_\epsilon(x)$
for a.e. $x \in \R^d$.  Also $w_\epsilon^\delta$ have the same
Lipschitz constant as $w_\epsilon$ and for any $\gamma > 0$, we
have
\begin{align}\label{St3s5}
\lambda w_\epsilon^\delta(x) &- \frac{\gamma^2}{2}
\ {\rm tr}\ (D^2 w_\epsilon^\delta(x)) - H(x, D w_\epsilon^\delta(x),
D^2 w_\epsilon^\delta(x)) \leq \rho_0(\epsilon)\nonumber\\
&+ g_\delta(x) + \frac{\gamma^2 d}{\epsilon},
\end{align}
where $g_\delta$ are uniformly continuous.  Let
\begin{equation*}
r_\delta(x, u_1, u_2) = r(x, u_1,u_2) + \rho_0(\epsilon)
+ \frac{\gamma^2 d}{\epsilon}.
\end{equation*}
We now assume $w(x) > \psi_2(x)$. Define the stopping time
$\hat{\tau}$ as before. Then $w_\epsilon^\delta(x) >
\psi_2(x)$ for sufficiently small $\epsilon$ and $\delta$.
Applying Ito's formula for $w_\epsilon^\delta$, we obtain for any
$T > 0$ and any stopping time $\theta$,
\begin{equation}\label{St3s6}
w_\epsilon^\delta(x) \leq E \left\lbrace \int_0^{T \wedge \theta \wedge
\hat{\tau}} {\rm e}^{-\lambda s} r_\delta(X^\gamma(s)) \,\, {\rm d}s
\!+\! {\rm e}^{-\lambda (T \wedge \theta \wedge \hat{\tau})}
w_\epsilon^\delta(X^\gamma(T \wedge \theta \wedge \hat{\tau}))
\right\rbrace,
\end{equation}
where $X^\gamma(t)$ is the solution of (\ref{St3s3}) with the
initial condition $X^\gamma(0) = x$.

By a standard martingale inequality, for any $\eta > 0$, we can
find a constant $R_\eta$  such that
\begin{equation*}
P \left(\sup_{0 \leq s \leq T} |X^\gamma(s)| \geq R_\eta \right)
\leq \eta.
\end{equation*}
Let $\Omega_\eta \subset \R^d$ be such that $|\Omega_\eta| \leq
\eta$ and $g_\eta \to 0$ uniformly on $B_{R_\eta} \setminus
\Omega_\eta$.  Using this we can find a local modulus $\rho_1$
such that
\begin{equation*}
\left\vert E \int_0^{T \wedge \theta \wedge \hat{\tau}} g_\delta
(X^\gamma(s)) \,\, {\rm d}s \right\vert \leq \rho_1(\delta, \gamma).
\end{equation*}
Using this in (\ref{St3s6}), we obtain
\begin{align}
w^\delta_\epsilon (x) &\leq E
\left\lbrace \int_0^{T \wedge \theta \wedge \hat{\tau}} {\rm e}^{-\lambda s}
r(X^\gamma(s))\,\, {\rm d}s + {\rm e}^{-\lambda (T \wedge \theta \wedge
\hat{\tau})} w^\delta_\epsilon(
X^\gamma(T \wedge \theta \wedge \hat{\tau})) \right\rbrace \nonumber\\
&\quad\ + \rho_2 (\delta, \gamma) + T \left( \rho_0(\epsilon) +
\frac{\gamma^2 d}{\epsilon} \right), \label{St3s7}
\end{align}
where $\rho_2$ is a local modulus. Now using moment estimates
\cite{S2}, we have
\begin{equation*}
E \left( \sup_{0 \leq s \leq T} |X(s) - X^\gamma(s)|^2 \right)
\leq C \gamma^2
\end{equation*}
for some constant $C >0$ which may depend on $T$ and $x$.  Now
using this in (\ref{St3s7}), and passing to the limits $\gamma \to
0$ and then $\delta, \epsilon \to 0$, we obtain
\begin{equation*}
w(x) \leq E \left\lbrace \int_0^{T \wedge \theta \wedge \hat{\tau}} {\rm
e}^{-\lambda s} r(X(s))\,\, {\rm d} s + {\rm e}^{-\lambda (T \wedge
\theta \wedge \hat{\tau})} w( X(T \wedge \theta \wedge \hat{\tau}))
\right\rbrace.
\end{equation*}
Now letting $T \to \infty$, as before, we obtain
\begin{equation*}
w(x) \leq V^-(x).
\end{equation*}
Similarly we can show that
\begin{equation*}
w(x) \geq V^+(x)
\end{equation*}
for all $x$ such that $\psi_1(x) > w(x)$. Thus the stochastic game
with stopping times has a value.

We now show that $(\hat{\theta}, \hat{\tau})$ is a saddle point
equilibrium. We need to prove this only when $\psi_2(x) < w(x) <
\psi_1(x)$. Using the foregoing arguments, we can show that
\begin{equation*}
w(x) = E \left[ \int_0^{\hat{\theta}\wedge \hat{\tau}} {\rm e}^{-\lambda s}
r(X(s)) \,\, {\rm d}s + {\rm e}^{-\lambda \hat{\theta} \wedge \hat{\tau}}
w(X(\hat{\theta}\wedge\hat{\tau})) \right] = R(x, \hat{\theta},
\hat{\tau}).
\end{equation*}
Clearly $(\hat{\theta}, \hat{\tau})$ constitutes a saddle point
equilibrium.\hfill \ab
\end{proof}

\begin{rmk}
{\rm The above result generalizes the optimal stopping time problem for
degenerate diffusions studied by Menaldi \cite{Me}.  Menaldi has
characterized the value function in the optimal stopping time
problem as a maximal solution of the corresponding variational
inequalities with one sided constraint.  He has used the
penalization arguments to obtain his results. Here we have used
the method of viscosity solutions to generalize the optimal
stopping time problem to stochastic games with stopping times.  We
also wish to mention that Stettner \cite{St} has studied
stochastic games with stopping times for a class of Feller Markov
processes.  He has employed a penalization argument using
semigroup theory.  Thus our approach is quite different from that
of Stettner.}
\end{rmk}

\section{Conclusions}

We have studied a class of  second order variational inequalities
with bilateral constraints. Under certain conditions, we have
showed the existence and uniqueness of viscosity solutions by
transforming the variational inequalities to HJI equations
corresponding to a stochastic differential game. Here we have
confined our attention to a particular form of $H$ which arises in
stochastic differential games of mixed type.  A general form of
$H$ can be reduced to this particular form by a suitable
representation formula as in \cite{Ka}. Thus our probabilistic
method can be used to prove the existence and uniqueness of
viscosity solutions for more general class of second order
nonlinear variational inequalities with bilateral constraints.

As an application, we have showed the existence of a value and a
saddle point equilibrium for a stochastic  game with stopping
times.  We now give a brief description of a stochastic
differential game of mixed type where each player uses both
continuous control and stopping times.

Consider the following controlled stochastic differential equation
\setcounter{equation}{0}
\begin{equation}\label{St1s1}
\hskip -2pc\left. \begin{array}{rcl}
{\rm d}X(t) & = & b(X(t), u_1(t), u_2(t)){\rm d}t
+ \sigma(X(t), u_1(t), u_2(t)){\rm d} W(t),\quad t>0 \\[.2pc]
X(0) & = & x
\end{array} \right\rbrace.
\end{equation}
Here $W(t)$ is a standard $d$-dimensional Brownian motion on
an underlying complete probability space $(\Omega, {\mathcal F},
P)$ and $\sigma$ is the non-negative square root of $a$.  The
processes $u_1(t)$ and $u_2(t)$ are $U_1$ and $U_2$-valued
control processes of players $1$ and $2$ respectively which are
progressively  measurable  with respect to the $\sigma$-field
generated by $W(t)$. Let $\lambda
> 0$ be the discount factor.  The payoff function is given by
\begin{align}
R(x, u_1(\cdot), \theta, u_2(\cdot), \tau) &= E \left[ \int_0^{\theta
\wedge \tau} {\rm e}^{-\lambda s}
r(X(s), u_1(s), u_2(s)) \,\, {\rm d}t \right. \nonumber \\
&~~~ + {\rm e}^{-\lambda (\theta \wedge \tau)} [ \psi_1(x(\theta))
\chi_{\theta < \tau} + \psi_2(x(\tau)) \chi_{ \tau \leq \theta}] \bigg\rbrack,
\label{St1s2}
\end{align}
where $\theta, \tau$ are the stopping times with respect to the
filtration generated by $W(\cdot)$.


An admissible control for player $1$ is a map $u_1(\cdot): [0,
\infty) \to U_1$ which is progressively measurable with respect to
the $\sigma$-field generated by $W(\cdot)$. The set of all
admissible controls for player $1$ is denoted by $\CA_1$.
Similarly an admissible control for player $2$ is defined. Let
$\CA_2$ denote the set of all admissible controls for player $2$.

We identify two  controls $u_1(\cdot), \tilde{u}_1(\cdot)$ in
$\CA_1$  on $[0, t]$ if $P(u_1(s) = \tilde{u}_1(s) \mbox{ for a.e.
} s \in [0, t]) =1$. Similarly we  identify the controls in
$\CA_2$.

An admissible strategy  for player $1$ is a map $\alpha: \CA_2 \to
\CA_1$ such that if $u_2 = \tilde{u}_2$ on $[0, s]$ then
$\alpha[u_2] = \alpha[\tilde{u}_2]$ for all $s\in [0, \infty)$.
The set of all admissible strategies for player $1$ is denoted by
$\Gamma$.  Similarly admissible strategies for player $2$ are
defined. Let $\Delta$ denote the set of all admissible strategies
for player $2$.

Let ${\mathcal S}$ denote the set of all  stopping times. Let
$\hat{\Gamma}$ denote the set of all  non-anticipating maps
$\hat{\alpha} : \CA_2  \to \CA_1 \times {\mathcal S}$ for player
$1$. Similarly let $\hat{\Delta}$ denote the set of all
non-anticipating maps  $\hat{\beta} : \CA_1 \to \CA_2 \times
{\mathcal S}$  for player $2$. Player $1$ tries to minimize $R(x,
u_1(\cdot), \theta, u_2(\cdot), \tau) $ over his admissible
control $u_{1}(\cdot)$ and stopping times $\theta$, whereas player
$2$ tries to maximize the same over his admissible control
$u_2(\cdot)$ and stopping times $\tau$. We now define the upper
and lower value functions of stochastic differential game of mixed
type. Let
\begin{align*}
V^+(x) &= \sup_{\hat{\beta} \in \hat{\Delta}} \inf_{u_1 \in \CA_1,
\theta \geq 0} R(x, u_1(\cdot), \theta, \hat{\beta}[u_1](\cdot)),\\
V^-(x) &= \inf_{\hat{\alpha} \in \hat{\Gamma}} \sup_{u_2 \in \CA_2,
\tau \geq 0} R(x, \hat{\alpha}[u_2](\cdot), u_2(\cdot), \tau).
\end{align*}
The functions $V^+$ and $V^-$ are respectively called upper and
lower value functions of the stochastic differential game of mixed
type. This differential game is said to have a value if both upper
and lower value functions coincide.  Now we make the following
conjecture.

\begin{conjecture} The value functions $V^+$
and $V^-$ are unique viscosity solutions of (\ref{St2s1}) and
(\ref{St2s2}) respectively in the class of bounded continuous
functions.

Note that the above conjecture is true for the special case
treated in \S3. Analogous results also holds when the matrix
$a$ is independent of the control variables and is uniformly
elliptic \cite{BL,Fr}.
\end{conjecture}

\section*{Acknowledgement}

This work is supported in part by IISc--DRDO Program on Advanced
Engineering Mathematics, and in part by NSF under the grants ECS-0218207
and ECS-0225448: also in part by Office of Naval Research through the
Electric Ship Research and Development Consortium. A part of this work
was done when this author was a CSIR Research Fellow at the Department
of Mathematics, Indian Institute of Science, Bangalore. The authors wish
to thank the anonymous referee for some corrections and improvements.

\end{document}